\documentclass[12pt,reqno]{amsart}

\usepackage[utf8]{inputenc}
\usepackage{amsmath} 
\usepackage{amsthm} 
\usepackage{amssymb} 
\usepackage{amscd} 

\usepackage{emptypage}
\usepackage{enumitem} 
\usepackage{mathtools} 
\usepackage{mathrsfs} 
\usepackage{upgreek}

\usepackage{hyperref,cleveref,autonum}

\usepackage{esint} 
\usepackage{graphicx,caption} 
\usepackage{float} 

\usepackage{pgf,tikz,pgfplots}
\usetikzlibrary{arrows}
\usetikzlibrary{intersections}
\usetikzlibrary{fadings}
\usepackage{wrapfig}

\usepackage[margin=0.9in]{geometry}
\newtheorem{theorem}{Theorem}[section]
\newtheorem*{theorem*}{Theorem}
\newtheorem{proposition}[theorem]{Proposition}

\theoremstyle{definition}

\theoremstyle{remark}
\newtheorem{remark}[theorem]{Remark}
\newtheorem*{remark*}{Remark}

\newcommand\beqc[1]{\left\{\begin{array}{#1}}
	\newcommand\eeqc{\end{array} \right.}

\def\bmatrix{\begin{pmatrix}}
	\def\ematrix{\end{pmatrix}}

\numberwithin{equation}{section}
\title[Conservation Laws and Energy Cascade]{New Conservation Laws and Energy Cascade for 1d Cubic NLS and the Schr\"odinger map}

\author[V. Banica]{Valeria Banica}
    \address[V. Banica]{$^{12}$ \newline $^1$  Sorbonne Universit\'e, CNRS, Universit\'e de Paris, Laboratoire Jacques-Louis Lions (LJLL), 4 place Jussieu, 75005 Paris, France, \newline $^2$Institut Universitaire de France (IUF).} 
\email{Valeria.Banica@sorbonne-universite.fr}

\author[L. Vega]{Luis Vega\\\\ \it{Dedicated to Carlos Kenig}}
\address[L. Vega]
	{$^{12}$ \newline $^1$ Basque Center for Applied Mathematics, Alameda Mazarredo 14, 48009 Bilbao, Spain, 
	\newline $^2$ Departamento  de  Matem\'aticas,  Universidad  del  Pais  Vasco/Euskal Herriko Unibertsitatea,  Aptdo.   644,  48080 Bilbao, Spain.
	}
\email{lvega@bcamath.org}

\date{\today}

\begin{document}

\begin{abstract}
We review some recent results concerning the Initial Value Problem of 1d-cubic non-linear Schr\"odinger equation (NLS) and other related systems as the Schr\"odinger Map. For the latter we prove the existence of a cascade of energy. Finally, some new examples of the Talbot effect at the critical level of regularity are given.  
\end{abstract}

\maketitle



\section{Introduction}
In these pages we review some recent results concerning the 1d-cubic non-linear Schr\"odinger equation (NLS) and other related systems. One of the main objectives is to explain in which sense
\begin{equation}
\label{monster}
u_M(x,t)=c_M\sum\limits_k e^{itk^2+ikx},
\end{equation}
for some constant $c_M$, is a ``solution'' of 1d--cubic NLS and to show the variety of phenomena it induces. Moreover, we will explain that it has a geometrical meaning due to its connection with  the Binormal Flow (BF) and the Schr\"odinger map (SM). 
Finally, we will explain how the so-called Talbot effect in Optics, that is mathematically described by $u_M$, is also present in the non-linear setting with data at the critical level of regularity. 

Altogether we are speaking about a family of PDE problems. Consider first NLS which is a complex scalar equation with a cubic non-linear potential:
\begin{equation}
\label{NLS}
\left\{\begin{split}\partial_t u&=i\left(\partial_x^2 u+(|u|^2-\mathcal{M}(t))\right)u,\qquad \mathcal{M}(t)\in\mathbb R,\\
u(x,0)&=u_0(x),\qquad x\in\mathbb R.
\end{split}\right.
\end{equation}

 Let us introduce next  SM. Calling $T(x,t)$ a unit vector in $\mathbb R^3$ the Schr\"odinger Map  onto the sphere is given by
\begin{equation}
\label{SM}\partial_tT=T\land \partial_x^2T.
\end{equation}

Finally, observe that the vector $T(x,t)$ can be seen at any given time as the tangent vector of a 3d-curve $\chi (x,t)$
$$\partial_x\chi= T,
$$
with $\chi$ a solution of 
\begin{equation}
\label{BF}\partial_t\chi=\partial_x\chi\land \partial_x^2\chi.
\end{equation}
Da Rios \cite{DaR} proposed \eqref{BF} as a simplified model that describes the evolution of vortex filaments. Remember that from Frenet equation $T_x=\kappa n$ with $\kappa$ the curvature of the curve and $n$ the normal vector. Hence
$$\partial_t\chi=\partial_x\chi\land \partial_x^2\chi=\kappa b$$
with $b$ denoting the binormal vector. That is the reason why sometimes the system of PDEs \eqref{BF} is called the Binormal Flow.

The connection of the two systems BF and SM with \eqref{NLS} was established by Hasimoto in \cite{Ha} thorugh a straightforward computation. This computation is slightly simplified if instead of the Frenet frame one uses the parallel one. This is given by vectors $(T,e_1,e_2)$ that satisfy

\begin{equation}
\label{parallel}
\begin{aligned}
T_x&=\hphantom{-\alpha T+}\alpha e_1+\beta e_2\\
e_{1x}&=-\alpha T\\
e_{2x}&=-\beta T.
\end{aligned}
\end{equation}
Defining
$$u=\alpha + i \beta$$
it is proved that $u$ solves \eqref{NLS} for some given $\mathcal{M}(t)$. 

 Observe that  $u$, the solution of \eqref{NLS}, gives the curvature and the torsion of $\chi$. More concretely,
$$|u|^2=\alpha^2+\beta^2=\kappa^2$$
and
\begin{equation}
\label{Hasimoto}
u(x,t)=\kappa e^{i\int_0^t \tau(r)\,dr},
\end{equation}
with $\bold\tau$ denoting the torsion. 
\eqref{Hasimoto} is usually called Hasimoto transformation.

A relevant simple example is 
\begin{equation}
\label{uo}
u_o=c_o\frac 1{\sqrt t} e^{i x^2/4t},\qquad \mathcal{M}(t)=\frac{c_o^2}{t},
\end{equation}
 which is related to the self-similar solutions of SM and BF. Formally $u_o(x,0)=c_o\delta$, and the corresponding $\chi$ has a corner at $(x,t)=(0,0).$ We will sometimes refer to this solution as either a fundamental brick or a coherent structure \cite{PL}. 
Our main interest is to consider rough initial data as polygonal lines and  regular polygons. As we will see, the latter are related to \eqref{monster} and therefore $u_M$ could be understood as a superposition of infinitely many simple solutions $u_o (x- j)$ centred at the integers. As a consequence the curve obtained from $u_M$ can be seen as an interaction of these coherent structures (see https://www.youtube.com/watch?v=fpBcwuY57FU).

It is important to stress that to obtain $\chi$ from $T$, besides integrating in the spatial variable the parallel system \eqref{parallel},  one has to find the trajectory in time followed by one point of $\chi_o$. This is not obvious even for \eqref{uo}, see \cite{GRV}. It turns out that to compute that trajectory of, say, one corner of a regular polygon is rather delicate, and the corresponding curve can be as complicated as those appearing as graphs of the so-called Riemann's non differentiable function:
\begin{equation}
\label{Riemann}\sum\limits_k\frac{e^{itk^2}-1}{k^2}.
\end{equation}
See  \cite{BVArma} for more details.

We will review some recent results regarding the IVP for \eqref{NLS} in section 2 and  section 3. In particular, we will show the existence of three new conservation laws \eqref{CL1}, \eqref{CL2}, and \eqref{CL3}, valid at the critical level of regularity. As it is well known \eqref{NLS}, and as a consequence also \eqref{SM} and \eqref{BF}, are completely integrable systems with infinitely many conservation laws that start at a subcritical level of regularity, $L^2$. For the others laws more regularity, measured in the Sobolev class, is needed. In section 4 we will recall some work done on the transfer of energy for the Schr\"odinger map \eqref{SM}. Finally, in section 5 we revisit the Talbot effect and modify some examples obtained in \cite{BV5} to establish a connection with some recent work on Rogue Waves given in \cite{GGKS}.

\section{The Initial Value Problem}

We start with the IVP associated to NLS equation \eqref{NLS}:
$$\left\{\begin{array}{cc}\partial_t u=i\left(\partial_x^2 u+(|u|^2-\mathcal{M}(t))\right)u,\qquad \mathcal{M}(t)\in\mathbb R,\\
u(x,0)=u_0(x),\qquad x\in\mathbb R.
\end{array}\right.$$
We are interested in initial data which are at the critical level of regularity. There are two symmetries that leave invariant the set of solutions that we want to consider. One is the scaling invariance: if $u$ is a solution of \eqref{NLS}, then
\begin{equation}
\label{scaling}\lambda>0\qquad u_\lambda (x,t)=\lambda u(\lambda x,\lambda^2t),
\end{equation}
is also a solution of \eqref{NLS} with $\lambda^2\mathcal{M}(t)$ instead of $\mathcal{M}(t)$.
The second one is the so-called Galilean invariance (i.e. translation invariance in phase space):
If 
\begin{equation}
\label{galilean}\nu\in\mathbb R\qquad u^{\nu} (x,t)=e^{-it\nu^2+i\nu x}u(x-2\nu t,t),
\end{equation}
then $u^{\nu}$ is also a solution of \eqref{NLS} with the same $\mathcal{M}(t).$
Hence, we want to work in a functional setting where the size of the initial data does not change under the scaling and Galilean transformations.

Let us review the classical results on NLS. The well-posedness of 1D cubic NLS on the full line and on the torus was firstly done  in \cite{T} and \cite{Bo} for data in $L^2$. Observe that the space $L^2(\Bbb R)$, although is invariant by Galilean symmetry \eqref{galilean}, misses the scaling \eqref{scaling} by a power of $1/2$ in the homogeneous Sobolev class $\dot H^s$. In fact, the critical exponet for scaling is $s=-1/2$ which is not invariant under the Galilean symmetry. The first result obtained beyond the $L^2(\Bbb R)$ theory was given in \cite{VV} using some spaces of tempered distributions built on the well known Strichartz estimates. Later on in \cite{G} well-posedness is studied in the Fourier-Lebesque spaces that we denote by $ \mathcal{F}L^p $. These are spaces where the Fourier transform is bounded in $ L^p(\Bbb R) $ and therefore leave invariant \eqref{galilean}. Moreover $ \mathcal{F}L^{\infty} $ is also scaling invariant and therefore critical according to our definition.    In \cite{GH} local well-posedness, also under periodic boundary conditions,  was shown in $ \mathcal{F}L^p $ with $ 2 < p < + \infty $.

In the setting of Sobolev spaces of non-homogeneous type the progress has been remarkable. On one hand, there is ill-posedness, in the sense that a data to solution map which is uniformly continuous does not exist in $ H^{s} $ with $ s < 0, $ and even some growth of the Sobolev norms has been proved,  \cite{KPV1},\cite{CCT}, \cite{CK},\cite{K},\cite{Oh}. On the other hand, it has been shown  in \cite{HKV} well-posedness in $ H^{s} $ for $ s > -1/2 $.   A weaker notion of continuity for the data to solution map is used.  

We will focus our attention in the symmetry of the translation invariant in phase (Fourier) space. We will consider solutions of \eqref{NLS} such that  
\begin{equation}
\label{omega}\omega(\xi, t):={e^{it\xi^2}\widehat u(\xi,t)\,\text{ is }\,2\pi\text{--periodic.}}
\end{equation}
Here
$\widehat u$ denotes the Fourier transform of $u$,
$$\int_\Bbb R e^{-ix\xi}u(x)\,dx.$$
To prove that this periodicity is preserved by the evolution is not completely obvious and it is a relevant property of  \eqref{NLS}. It can be proved writing the equation for $\omega$, $\mathcal{M}(t)=0)$ in \eqref{NLS}:
$$
\partial_t \omega(\eta, t)=\frac{i}{8\pi^3}e^{-it\eta^2}\int\int_{\xi_1+\xi_2+\xi_3-\eta=0} e^{it(\xi_1^2-\xi_2^2+\xi_3^2)} \omega(\xi_1)\bar\omega(\xi_2)\omega(\xi_3) \,d\xi_1d\xi_2d\xi_3.
$$
Under the condition $
\xi_1-\xi_2+\xi_3-\eta=0$, we get
$$\xi_1^2-\xi_2^2+\xi_3^2-\eta^2=2(\xi_1-\xi_2)(\xi_1-\eta).
$$
The last quantity is invariant under translations so that the periodicity is formally preserved. Interestingly this calculation does not work for general dispersive systems as for example for modified KdV.

One of the three new conservation laws is precisely
\begin{equation}
\label{CL1}\int_0^{2\pi}|\omega(\xi, t)|^2 \, d\xi= constant.
\end{equation}
This can be seen writing $\omega(\xi,t)=\sum\limits_j A_j(t) e^{ij\xi}$ and looking for the ODE system that the Fourier coefficients $A_j$ have to satisfy. 
Historically, our approach to this question has been different and this is what we explain next.

Following Kita in \cite{Kita}, we considered the  ansatz
\begin{equation}
\label{ansatz}
u(x,t)=\sum\limits_j A_j(t) e^{it\partial_x^2}\delta(x-j),
\end{equation}
and therefore
$$\widehat u(\xi,t)=e^{-it\xi^2}\sum\limits_j A_j(t) e^{ij\xi}.$$
If we define
\begin{equation}
\label{V1}
V(y,\tau)=\sum\limits_j B_j(\tau)e^{ ijy},
\end{equation}
then
$$\begin{aligned}
u(x,t)&=\displaystyle \frac{1}{(it)^{1/2}}\sum\limits_j A_j(t)e^{i\frac{(x-j)^2}{4t}}\\
&=\displaystyle \frac{1}{(it)^{1/2}}e^{i\frac{|x|^2}{4t}}\sum\limits_j A_j(t)e^{i\frac{j^2}{4t}-i\frac{x}{2t}j}\\
&:{=\displaystyle \frac{1}{(it)^{1/2}}e^{i\frac{|x|^2}{4t}}\overline V\left(\frac{x}{2 t},\frac 1t\right),
}\end{aligned}$$
with
\begin{equation}
\label{B}
B_j(\tau)=\overline{A_j}\left(\frac 1\tau\right)e^{-i\frac \tau 4 j^2}.
\end{equation}
Finally, doing the change of variables
$$y=\displaystyle\frac{x}{2 t,}\qquad \tau=1/t,
$$ 
we easily obtain that $V$ solves
\begin{equation}
\label{V}
\partial_{\tau}V=i\left(\partial_y^2+\displaystyle\frac 1\tau(|V|^2-m)\right)V\qquad;\qquad \displaystyle m(\tau)=\frac 1\tau \mathcal{M}(\frac 1{\tau}).
\end{equation}
We actually have that $V$ is a pseudo-conformal transformation of $u$.
\begin{remark}
\label{rem2}
\begin{enumerate}

\item Observe that formally solutions of \eqref{V} remain periodic if they are periodic at a given time. That means that given the Fourier coefficients $B_j$ and using \eqref{B} to define $A_j$ we conclude that the periodicity of $\omega(\xi, t)=e^{it\xi^2}\widehat u(\xi,t)$ is also formally preserved. 
\item There is a singularity at $\,t=0\,$ artificially created for the change of variable. Hence a very natural question is if $u$ can be continued for $t\le 0$. This issue becomes a question about the scattering of the solutions of \eqref{V}.
\item If $V(1)=c_o,$ and $m=c_o^2$ then $V(\xi,t)=c_o$  for all t. The corresponding solution is the fundamental brick \eqref{uo}
$$
u(x,t)=u_o=\displaystyle\frac {c_o}{\sqrt t} e^{i\frac{|x|^2}{4t}},
$$
and $u_0=c_o\delta$. This implies that unless we include the term $\mathcal{M}(t)=c_o^2/t$ in \eqref{NLS} the IVP for the Dirac delta is ill--posed, something observed in \cite{KPV1}. It was proved in \cite{BV2009}  and \cite{BV2015} that even if this term is added and one looks for solutions of the type $V=c_o+z\,$ with $\,z\,$ small with respect to $\,a\,$ the corresponding $u$ of \eqref{NLS} cannot be defined for $t=0$. 
\item From \eqref{V} it immediately follows that 
\begin{equation}
\label{CL1bis}m_0=\int_{0}^{2\pi} |V(\xi,\tau)|^2 d\xi=\sum\limits_{j}|B_j(\tau)|^2
\end{equation} is formally constant for $\tau>0$. And from \eqref{B} we also get that $\sum\limits_j|A_j(t)|^2=\sum\limits_j|B_j(1/t)|^2=\int_0^{2\pi}|\omega(\xi,t)|^2\,d\xi$ remains constant, which is \eqref{CL1}.
\end{enumerate}
\end{remark}

\section {Conservation laws}


In \cite{BV5} a first result on the IVP \eqref{NLS} within the functional setting we have just described was obtained with the ansatz \eqref{ansatz}. In fact, the solution $u$ is written as
\begin{equation}
\label{ansatz2} 
u(x,t)=\sum\limits_j A_j(t) e^{i\Phi_j(t)}e^{it\delta_x^2}\delta(x-j)\qquad \phi_j(t)=e^{i\frac{|a_j|^2}{8\pi}\log t}.
\end{equation}
If for any $a_j$ we write 
\begin{equation}
\label{R_j}A_j(t)=a_j+R_j(t),
\end{equation}
then an infinite ODE system for the $R_j$'s can be easily obtained. The corresponding solution is constructed through a fixed point argument in an appropriately chosen space which among other things implies that
$$ R_j(0)=0.$$
The condition on the data is that  $\sum_j |a_j|$
is finite (i.e. $a_j\in l^1$) but not necessarily small. The result is local in time. A global result is obtained by assuming the extra condition 
\begin{equation}
\label{2moment}\,\sum\limits_j j^2|a_j|^2<+\infty,
\end{equation}
whose evolution $\sum_j |A_j(t)|^2$ is easy to determine as we explain next.

First of all, it is much more convenient to work with $V$ defined in \eqref{V1}, solution of \eqref{V}. Then, it is easy to compute the ODE system that the Fourier coefficients $B_j$'s of $V$ have to satisfy:
\begin{equation}\label{Bjode}i\partial_\tau B_k(\tau)=\frac{1}{ \tau}\sum_{k-j_1+j_2-j_3=0}e^{-i\tau(k^2-j_1^2+j_2^2-j_3^2)}B_{j_1}(\tau)\overline{B_{j_2}(\tau)}B_{j_3}(\tau).
\end{equation}
Observe that the condition $ k-j_1+j_2-j_3=0$ has to be satisfied. Calling 
\begin{equation}
\label{w}w_{k,j_1,j_2}:=k^2-j_1^2+j_2^2-j_3^2=2(k-j_1)(j_1-j_2), 
\end{equation}
the resonant set is given by  $w_{k,j_1,j_2}=0.$ Thus $c(k)-c(j_1)+c(j_2)-c(j_3)$ vanishes on the resonant set for any real function $c$. 
Let us introduce the non-resonant set
$$NR_k=\{(j_1,j_2,j_3),\, k-j_1+j_2-j_3=0,\, k^2-j_1^2+j_2^2-j_3^2\neq 0\}.$$
Then, for any such a function $c$ we have
\begin{equation}\label{Bjcons}
\begin{split}\frac{d}{d\tau} \sum_k c(k)|B_k(\tau)|^2&=\frac{1}{2 i\tau }\sum_{k-j_1+j_2-j_3=0}(c(k)-c(j_1)+c(j_2)-c(j_3))e^{-i\tau w_{k,j_1,j_2}}B_{j_1}(\tau)\overline{B_{j_2}(\tau)}B_{j_3}(\tau)\overline{B_k(\tau)}\\
&=\frac{1}{2i\tau}\sum_{k;NR_k}(c(k)-c(j_1)+c(j_2)-c(j_3))e^{-i\tau w_{k,j_1,j_2}}B_{j_1}(\tau)\overline{B_{j_2}(\tau)}B_{j_3}(\tau )\overline{B_k(\tau)}.
\end{split}
\end{equation}
Relevant examples  are $c(j)=1$ that gives \eqref{CL1} ,the $L^2$ conservation law already mentioned, and $c(j)=j$, that yields a second conservation law:
\begin{equation}
\label{CL2}
\sum\limits_j j|B_j(\tau)|^2=\sum\limits_j j|A_j(1/\tau)|^2\text{\,\,is constant.}
\end{equation}
The final example is $c(j)=j^2$, cf. \eqref{2moment}. In this case the corresponding quantity does not remain constant and its derivative is better understood in terms of $V$. Calling
$$E(\tau)=\int |\partial_y V (y,\tau)|^2-\frac1{4\tau}(|V(y,\tau)|^2-m)^2\,dy,
$$
we get
$$\frac d{d\tau}E(\tau)=\frac1{4\tau^2}\int(|V(y,\tau)|^2-m)^2\,dy.
$$
Finally observe denoting $m_0=\sum\limits _j |B_j|^2$ equation \eqref{Bjode} can be written as
\begin{equation}\label{Bjodev}i\frac d{d\tau} B_k(t)=\frac{1}{ \tau}\sum_{NR_k}e^{-itw_{k,j_1,j_2}}B_{j_1}(\tau)\overline{B_{j_2}(\tau)}B_{j_3}+\frac 1\tau\left(2m_0-|B_k(t)|^2\right)B_k(\tau).\end{equation}

The next step about the IVP \eqref{NLS}-\eqref{ansatz2} was given in \cite{BVArma}, where the Picard iteration is done measuring more carefully the first iterate. Particular attention is given to the example
\begin{equation}
\label{example1}
a_j=1\,\,\text{for}\,\,|j|\le N\quad \text{and zero otherwise},
\end{equation}
see \cite{BVArma}.

Finally in \cite{BrV}, Bourgain's approach \cite{Bo} is followed. This amounts to use the Sobolev spaces of the coefficients $B_j(\tau)$'s. The results in that paper can be summarized as follows, for initial datum in $ l^p$, $  p \in (1,+\infty) $: 

\begin{enumerate}
	
	\item  Local well-posedness with a smallness assumption in $ l^p $ for the initial datum: for any $ T > 0 $, there exists $ \epsilon(T) > 0 $ such that if the $ l^p $ norm of the initial datum $ \{a_{j}\} $ is smaller then $ \epsilon(T) $, then there exists a unique solution of \eqref{ansatz2} in $ [0, T] $ in an appropriate sense. 
	
\item  Local well-posedness with a smallness assumption in $ l^{\infty} $ for the initial datum: if the $ l^{\infty} $ norm of $\{ a_j \}$ is small enough then there exists a  time $ T(\|\alpha\|_{\l^{\infty}},\|\alpha\|_{\l^p}) $ such that a unique solution of \eqref{ansatz2} exists in $ [0,T] $ in an appropriate sense.  
\item    For $ p = 2 $, global in time well-posedness with a smallness assumption in $ l^{\infty} $ for the initial datum. As it can be expected this result follows from (2) and the $ l^2 $ conservation law. The smallness condition comes from the linear term that is treated as a perturbation. We don't know if this smallness condition can be removed.  
\end{enumerate}

For establishing the third conservation law we have to observe that $w$ given in \eqref{w} is invariant under translations. This implies that if $B_{j+M}=B_j$ at a given time the property is formally preserved for all time and therefore
\begin{equation}
\label{CL3}m_0=\sum_{j=1}^M|a_j|^2=\sum_{j=1}^M|A_j(t)|^2=\sum_{j=1}^M|B_j(1/t)|^2=constant.
\end{equation}
This conservation law is much stronger than \eqref{CL1} because it just assumes an $l^\infty$ condition on the $a_j$'s.

As a consequence, in \cite{BrV} an explicit solution of \eqref{Bjodev} for the relevant case
\begin{equation}
\label{constant}
c_M=a_j=B_j \qquad\text{for all}\ j
\end{equation}
is constructed as \eqref{monster}.
\newline\newline

From the results reviewed in this section and the ansatz \eqref{ansatz2} we conclude that the IVP for $\,u\,$ at $\,t=0\,$ is ill--posed due to a loss information of the phase. As it was proved in \cite {BV5} this loss is irrelevant when \eqref {NLS} is understood in connection to BF and SM. For example, for BF the solutions can be perfectly defined at $t=0$ as a polygonal line that, except for \eqref{constant}, tend to two straight lines at infinite. Moreover, the behavior close to a corner is determined by a self-similar solution \eqref{uo}. This self-similar solution, and the precise theorem given about them in \cite{GRV}, gives the necessary information at $t=0$ so that the flow can be continued for $t<0$. A crucial ingredient in this process is the precise relation established in \cite{GRV} between $c_o$ and the angle $\theta_o$ of the corresponding corner, namely
\begin{equation}
\label{theta}
\sin\theta_o=e^{-\frac{\pi c_0^2}{2}}
\end{equation}

For a regular polygon with $M$ sides the angle is $\theta_M=2\pi/M$. Chossing $c_M$  in such a way that
$$\sin(\frac{2\pi}{M})=e^{-\frac{\pi c_M^2}{2}}$$
and using it in \eqref{constant}, we obtain a solution for the case of a regular polygon at the level of NLS. This choice is the one conjectured in \cite{DHV} based on the numerical simulations done in \cite{JeSm2}.

\section {Transfer of energy}

In section 3 we have constructed solutions of \eqref{NLS} whose energy density is well described in terms of $|\widehat u|^2$ as
$|\widehat u|^2=|\omega|^2$
with $\omega$ given in \eqref{omega}. Also remember that $\omega$ is related to $V$ through  \eqref{B} and \eqref{V}. This raises the question about up to what extent $\big|\widehat{T_x}(\xi,t)\big|^2 d\xi$ can be considered a density energy. A hint that suggests a positive answer is given in \cite{BVCPDE} where the following identity is proved
$$\displaystyle\int_0^{2\pi} |V(\xi,t)|^2 d\xi=\lim\limits_{n\to\infty}\int_{2\pi n}^{2\pi(n+1)}
\big|\widehat{T_x}(\xi,t)\big|^2 d\xi.$$

Even though it was proved in section 2 and section 3 that
$$\displaystyle\int_0^{2\pi} |V(\xi,t)|^2 d\xi= constant$$
and therefore, that there is no flux  of energy for $\widehat u$ at least for $0<t<1$ the situation for $T$ is different. It was proved in  \cite{BVIHP} that there is some cascade of energy. More concretely we have the following result that was motivated by some numerical experiments done in \cite{DHV2}.

\begin{theorem}
 Assume

$$\left\{\begin{array}{cl}
a_{-1}=a_{+1}\ne 0,\\
a_j=0&\text{otherwise.}
\end{array}\right.$$
 Then there exists $\,c>0$ such that
$$\sup\limits_\xi\big|\widehat {T_x} (\xi,t)\big|^2\ge  \sup\limits_{\xi\in B(\pm\frac 1t,\sqrt{t})}\big|\widehat {T_x} (\xi,t)\big|^2\ge c|\log t|\qquad t>0.$$
\end{theorem}

This type of energy cascade is an alternative to the ones in \cite{Iteam} and \cite{HPTV}. 

Recall that if $u=\alpha + i \beta$ is the solution of \eqref{NLS}  then $T$ can be obtained from \eqref{parallel}.  A simple calculation which can be found for example in \cite {BVCPDE} gives 
\begin{equation}
\label{Tt}
\begin{aligned}
T_t&=\hphantom{-\alpha_x T\qquad\qquad\,\,\,}-\beta_x e_1+\hphantom{\qquad\qquad\,\,.}\alpha_x e_2\\
e_{1_t}&=-\alpha_x T\hphantom{(|u|^2-\mathcal{M}(t))e_{1}-\,}+((\alpha^2+\beta^2)-\mathcal{M}(t))e_2\\
e_{2_t}&=-\beta_x T-((\alpha^2+\beta^2)-\mathcal{M}(t))e_1.
\end{aligned}
\end{equation}

Notice that this is just a linear system of equations which is hamiltonian and that satisfies the three conservation laws \eqref{CL1}, \eqref{CL2}, and \eqref{CL3} given in the previous sections for $u$ of type \eqref{ansatz2}. Nevertheless Theorem 4.1 applies and therefore this system also exhibits a cascade of energy. 

\section{Talbot effect and Rogue Waves}

In this section we want to revisit the examples on the Talbot effect showed in \cite{BV5}. The Talbot effect is very well described by \eqref{monster}. As it will be shown below in \eqref{perevol} and \eqref{Talbotper}, the values of \eqref{monster} at times which are rational multiples of the period can be written in a closed formed: if the rational is $p/q$ then Dirac deltas appear at all the rationals of $\mathbb Z/q$, and the amplitudes are given by a corresponding Gauss sum. Going either backward or forward in time this gives a phenomenon of constructive/destructive interference that we think it is similar to the one exhibited in \cite{GGKS} related to the so-called Rogue Waves.

The example we propose is very similar to \eqref{example1}. Recall that the construction we do is perturbative  and therefore it always implies some smallness condition. This condition is measured in terms of $\sum |a_j|$ that can be small without the corresponding solution $u$ being small. For example, from \eqref{ansatz} it is immediate that at least for small times  the $L^\infty$ norm of $u$ is not small. Something similar can be said for the $L^1_{loc}$ norm. At this respect it is relevant to notice the definition of $u$ in terms of $V$ given in \eqref{V}. Observe that 
$$|u(x,t)|=\frac1{\sqrt {t}}\, |V(x/2t,1/t)|,$$
and therefore $L^1_{loc}$ grows with $t$.

We have the following result.

\begin{theorem}\label{proptalbotnl}{\bf{(Appearance of rogue waves)}} Let $0<\eta<\frac 14$ and let $p\in\mathbb N$ large. There exists $u_0$ with $\widehat{u_0}$ a $2\pi-$periodic function, located modulo $2\pi$ in $[-\eta \frac {2\pi}{ p},\eta \frac {2\pi}{ p}]$, such that the solution $u(t,x)$ of \eqref{NLS} obtained from $a_k=\widehat{u_0}(k)$ in \cite{BV5} satisfies the following property. For times $t_{p,q}=\frac 1{2\pi}\frac pq$ and $t_{\tilde p,\tilde q}=\frac 1{2\pi}\frac {\tilde p}{\tilde q}$ both of size $\frac 1{2\pi}$, with rational representation of type $q\approx p$, and $\tilde p\approx \tilde q\approx 1$ with $\tilde p<\tilde q$, $q,\tilde q$ odd numbers, on the interval $[-\frac 1{2\tilde q},\frac 1{2\tilde q}]$ we observe at time $t_{p,q}$ almost-periodic small waves while at time $t_{\tilde p,\tilde q}$ a localized large-amplitude structure emerges.

\end{theorem}

We start with a computation for the linear Schrödinger equation on the line, concerning the Talbot effect related to \eqref{monster}. 
\begin{proposition}\label{thTalbot}{\bf{(Talbot effect for linear evolutions)}} 
Let $0<\eta<\frac 14$, $p\in\mathbb N$ and $u_0$ such that $\widehat{u_0}$ is a $2\pi-$periodic function, located modulo $2\pi$  in $[-\eta \frac {2\pi}{ p},\eta \frac {2\pi}{ p}]$. For all $t_{p,q}=\frac 1{2\pi}\frac pq$ with $q$ odd and for all $x\in\mathbb R$ we define
$$\xi_{x} :=\frac {\pi q}p \,d(x,\frac1q\mathbb Z) \in[0,\frac \pi p).$$
Then, there exists $\theta_{x,p,q}\in\mathbb R$ such that
\begin{equation}\label{repr2}e^{it_{p,q}\Delta}u_0(x)=\frac {1}{\sqrt{q}} \,\widehat{u_0}(\xi_x)\, e^{-it_{p,q}\,\xi_x^2+ix\,\xi_x+i\theta_{x,p,q}}.\end{equation}
In particular $|e^{it_{p,q}\Delta}u_0|$ is $\frac 1q$-periodic and if $d(x,\frac1q\mathbb Z)>\frac {2\eta}{q}$ then $e^{it_{p,q}\Delta}u_0(x)$ vanishes. 
\end{proposition}

\begin{proof}
We start by recalling the Poisson summation formula $\sum _{k\in\mathbb Z}f(k)=\sum _{k\in\mathbb Z}\hat f(2\pi k)$ for the Dirac comb:
$$
(\sum _{k\in\mathbb Z}\delta_k)(x)=\sum _{k\in\mathbb Z}\delta(x-k)=\sum_{k\in\mathbb Z}e^{i2\pi kx},$$
as 
$$\widehat{\delta(x-\cdot)}(2\pi k)=\int_{-\infty}^\infty e^{-i2\pi ky}\delta(x-y)\,dy=e^{-i2\pi kx}.$$

The computation of the free evolution with Dirac comb data is
\begin{equation}\label{perevol}e^{it\Delta}(\sum _{k\in\mathbb Z}\delta_k)(x)=\sum_{k\in\mathbb Z}e^{-it(2\pi k)^2+i2\pi kx}.\end{equation}
For $t=\frac 1{2\pi}\frac pq$ we have (choosing $M=2\pi$ in formulas (37) combined with (42) from \cite{DHV})
\begin{equation}\label{Talbotper}e^{it\Delta}(\sum _{k\in\mathbb Z}\delta_k)(x)=\frac 1{q}\sum_{l\in\mathbb Z}\sum_{m=0}^{q-1}G(-p,m,q)\delta(x-l-\frac mq),\end{equation}
which describes the linear Talbot effect in the periodic setting. Here $G(-p,m,q)$ stands for the Gauss sum
$$G(-p,m,q)=\sum_{l=0}^{q-1}e^{2\pi i\frac{-pl^2+ml}{q}}.$$\\

Now we compute the free evolution of data $u_0$ with $\widehat{u_0}$ a $2\pi-$periodic function, i.e. $\widehat{u_0}(\xi)=\sum_{k\in\mathbb Z}\alpha_ke^{-ik\xi}$ and $u_0=\sum _{k\in\mathbb Z}\alpha_k\delta_k$:
$$e^{it\Delta}u_0(x)=\frac 1{2\pi}\int_{-\infty}^\infty e^{ix\xi}e^{-it\xi^2}\widehat{u_0}(\xi)\,d\xi=\frac 1{2\pi}\sum_{k\in\mathbb Z}\int_{2\pi k}^{2\pi(k+1)} e^{ix\xi-it\xi^2}\widehat{u_0}(\xi)\,d\xi$$
$$=\frac 1{2\pi}\int_0^{2\pi}\widehat{u_0}(\xi)\sum_{k\in\mathbb Z}e^{ix(2\pi k+\xi)-it(2\pi k+\xi)^2}\,d\xi=\frac 1{2\pi}\int_0^{2\pi}\widehat{u_0}(\xi)e^{-it\xi^2+ix\xi}\sum_{k\in\mathbb Z}e^{-it\,(2\pi k)^2+i2\pi k (x- 2t\xi )}\,d\xi.$$
Therefore,  for $t_{p,q}=\frac 1{2\pi}\frac pq$ we get using \eqref{perevol}-\eqref{Talbotper}:
$$e^{it_{p,q}\Delta}u_0(x)=\frac 1{q}\int_0^{2\pi}\widehat{u_0}(\xi)e^{-it_{p,q}\xi^2+ix\xi}\sum_{l\in\mathbb Z}\sum_{m=0}^{q-1}G(-p,m,q)\delta(x-2t_{p,q}\xi-l-\frac mq)\,d\xi.$$
For $q$ odd  $ G(-p,m,q)=\sqrt{q}e^{i\theta_{m,p,q}}$ for some $\theta_{m,p,q}\in\mathbb R$ so we get for $t_{p,q}=\frac 1{2\pi}\frac pq$
$$e^{it_{p,q}\Delta}u_0(x)=\frac {1}{\sqrt{q}}\int_0^{2\pi}\widehat{u_0}(\xi)e^{-it_{p,q}\xi^2+ix\xi}\sum_{l\in\mathbb Z}\sum_{m=0}^{q-1}e^{i\theta_{m,p,q}}\delta(x-2t_{p,q}\,\xi-l-\frac mq)\,d\xi.$$
For a given $x\in\mathbb R$ there exists a unique  $l_x\in\mathbb Z$ and a unique $0\leq m_x<q$ such that 
$$x-l_x-\frac {m_x}q \in[0,\frac 1q),\quad \xi_{x} =\frac {\pi q}p(x-l_x-\frac {m_x}q)=\frac {\pi q}p\,d(x,\frac1q\mathbb Z) \in[0,\frac \pi p).$$
We note that for $0\leq \xi<\eta \frac{2\pi}p$ we have $0\leq 2t\xi<\frac {1}{2q}$. 
As $\widehat u_0$ is located modulo $2\pi$ only in a neighborhood of zero of radius less than $\frac \eta {2\pi}{ p}$ then we get the expression \eqref{repr2}.

\end{proof}

\begin{proof}
{\it (of Theorem \ref{proptalbotnl})} We shall construct sequences $\{\alpha_k\}$ such that $\sum _{k\in\mathbb Z} \alpha_k \delta_k$ concentrates in the Fourier variable near the integers. To this purpose we consider, for $s>\frac 12$, a positive bounded function $\psi\in H^{s}$ with support in $[-1,1]$ and maximum at $\psi(0)=1$. We define the $2\pi$-periodic function satisfying
$$f(\xi):=p^{\beta}\psi(\frac p{2\pi \eta} \xi),\quad \forall \xi\in[-\pi,\pi],$$
with $\beta<\frac 12-\frac 32s$,
 introduce its Fourier coefficients:
$$f(\xi):=\sum _{k\in\mathbb Z} \alpha_k e^{ik\xi},$$
and the function 
$$u_0:=\sum_{k\in\mathbb Z} \alpha_k \delta_k.$$
In particular, on $[-\pi,\pi]$, we have $\widehat{u_0}=f$ and the restriction of $\widehat{u_0}$ to $[-\pi,\pi]$ has support included in a neighborhood of zero of radius less than $\eta\frac {2\pi} p$. We then get from \eqref{repr2}:
\begin{equation}\label{repr3}e^{it_{p,q}\Delta}u_0(x)=\frac {1}{\sqrt{q}} \,\widehat{u_0}(\frac {\pi q}p \,d(x,\frac1q\mathbb Z))\, e^{-it_{p,q}\,\xi_x^2+ix\,\xi_x+i\theta_{m_x}},\end{equation}
that
\begin{equation}\label{estpointw}|e^{it_{p,q}\Delta}u_0(0)|=\frac {1}{\sqrt{q}} \,|f(0)|=\frac {1}{\sqrt{q}} p^{\beta}\psi(0)=\frac {p^{\beta}}{\sqrt{q}},\end{equation}
\begin{equation}\label{estmod}\|e^{it_{p,q}\Delta}u_0\|_{L^\infty}\leq \frac {1}{\sqrt{q}} \,\|f\|_{L^\infty}=\frac {p^{\beta}}{\sqrt{q}},\end{equation}
and
\begin{equation}\label{estpointwfar}e^{it_{p,q}\Delta}u_0(x)=0, \mbox { if }d(x,\frac1q\mathbb Z)>\frac {2\eta}{q}.\end{equation}

We note that
$$\|\alpha_k\|_{l^{2,r}}^2=\sum_k|k|^{2r}|\alpha_k|^2=\|f\|_{\dot H^r}^2=\frac{p^{2(\beta+r-\frac 12)}}{(2\pi\eta)^{2(r-\frac 12)}}\|\psi\|_{\dot H^r}.$$
Since $\beta<\frac 12-s$ and $p$ is large, it follows that $\|\alpha_k\|_{l^{2,s}}$ is small enough so that we can use the results  in \cite{BV5} to construct a solution up to time  $t=\frac 1{2\pi}$ for \eqref{NLS} 
of type
$$u(t,x)=\sum_{k\in\mathbb Z}e^{i(|\alpha_k|^2-2\sum_j|\alpha_j|^2)\log t}(\alpha_k+R_k(t))e^{it\Delta}\delta_k(x).$$
Hence
$$\left|u(t_{},x)-e^{it_{}\Delta}(\sum_{k\in\mathbb Z}\alpha_k\delta_k)(x)\right|$$
$$\leq \sum_{k\in\mathbb Z}(1-e^{i(|\alpha_k|^2-2\sum_j|\alpha_j|^2)\log t_{}})\alpha_k e^{it_{}\Delta}\delta_k(x)+\sum_{k\in\mathbb Z}e^{i(|\alpha_k|^2-2\sum_j|\alpha_j|^2)\log t_{}}R_k(t_{})e^{it_{}\Delta}\delta_k(x)$$
$$\leq \frac{C}{\sqrt{t_{}}}\|\alpha_k\|_{l^{2}}^2\|\alpha_k\|_{l^{2,s}} +\frac{Ct^\gamma}{\sqrt{t}}\|\alpha_k\|_{l^{2,s}}^3=C(\eta)\frac{p^{3(\beta+s-\frac 12)}}{\sqrt{t}}(p^{-s}+t^\gamma),$$
for $\gamma<1$.

Therefore, in view of \eqref{estpointw},\eqref{estmod} and \eqref{estpointwfar} we have for times $t_{p,q}$ and $t_{\tilde p,\tilde q}$ both of size $\frac 1{2\pi}$, but with rational representation of type $q\approx p$ which is fixed to be large, and $\tilde p\approx \tilde q\approx 1$ with $\tilde p<\tilde q$, that:
\begin{itemize}
\item at time $t_{p,q}$ the modulus $|u(t_{p,q},x)|$ is a $\frac 1q$-periodic function of maximal amplitude $\frac 1{p^{\frac 12-\beta}}$ plus a remainder term of size $\frac 1{p^{3(-\beta-s+\frac 12)}}$, that is negligible provided that $\beta<\frac 12-\frac 32s$. So modulo negligible terms $|u(t_{p,q},x)|$ has plenty of $\frac 1p$-period waves of small amplitude $\frac 1{p^{\frac 12-\beta}}$,

\item at time $t_{\tilde p,\tilde q}$ the modulus $|u(t_{\tilde p,\tilde q},x)|$ is a $\frac 1{\tilde q}$-periodic function of maximal amplitude $\frac 1{\tilde p^{\frac 12-\beta}}$ plus a remainder term of size $\frac 1{\tilde p^{3(-\beta-s+\frac 12)}}$, that is again negligible provided that $\beta<\frac 12-\frac 32s$. So modulo negligible terms $|u(t_{\tilde p,\tilde q},x)|$ has in the interval $I:=[-\frac 1{2\tilde q},\frac 1{2\tilde q}]$ a wave of amplitude $\frac 1{\tilde p^{\frac 12-\beta}}$, and is upper-bounded by a smaller value on $I\setminus [-\frac {2\eta}{\tilde q},\frac {2\eta}{\tilde q}]$.

\end{itemize}
Therefore, observing what happens in the interval $I$ we have at time $t_{p,q}$ small almost-periodic waves while at time $t_{\tilde p,\tilde q}$ a localized large-amplitude (with respect to $\eta$) structure emerges.

\end{proof}

\begin{remark}
\begin{enumerate}
\item In the above argument we need $\eta$ to be small. As a consequence,  the $L^\infty$ norm and therefore the $L^1_{loc}$ norm of the  solution is small. This can be avoided by considering $u_{\lambda}=\frac1\lambda u(x/\lambda, t/\lambda^2$), where $u$ is any of the solutions constructed above. If $\lambda>1$ the $L^\infty$ norm grows, while for $\lambda<1$ the $L^1_{loc}$ norm around  the corresponding bump grows.
\item The size of the error can be made smaller following the ideas developed in \cite{BVArma}. This is due to the type of data  \eqref{example1} we are using. In this case the size of the first Picard iterate is indeed much smaller than the $l^1$ norm we are using in the above argument.
\end{enumerate}
\end{remark}

{\bf{Acknowledgements:}}  This research is partially supported by the Institut Universitaire
de France, by the French ANR project SingFlows, by ERCEA Advanced Grant 2014 669689
- HADE, by MEIC (Spain) projects Severo Ochoa SEV-2017-0718, and PGC2018-1228
094522-B-I00, and by Eusko Jaurlaritza project IT1247-19 and BERC program.

\end{document}